\documentclass[12pt,a4paper]{article}

\usepackage[dvips]{color}
\usepackage{color}
\usepackage{xcolor}
\usepackage{mathrsfs}
\usepackage{amsfonts}
\usepackage{mathtools}
\usepackage{stmaryrd}
\usepackage{amsmath}
\usepackage{amsthm}
\usepackage{amssymb}
\usepackage{bbm}
\usepackage{graphicx}
\usepackage{enumerate}

\makeatletter
\def\tank#1{\protected@xdef\@thanks{\@thanks
        \protect\footnotetext[0]{#1}}}
\def\bigfoot{

    \@footnotetext}
\makeatother

\newcommand{\DEQS}{\begin{eqnarray*}}
\newcommand{\EEQS}{\end{eqnarray*}}
\newcommand{\DEQZ}{\begin{eqnarray}}
\newcommand{\EEQZ}{\end{eqnarray}}

\theoremstyle{plain}

\newtheorem{thm}{Theorem}[section]
\newtheorem{lem}{Lemma}[section]
\newtheorem{prop}{Proposition}[section]

\newtheoremstyle{boldremark}
    {\dimexpr\topsep/2\relax} 
    {\dimexpr\topsep/2\relax} 
    {}          
    {}          
    {\bfseries} 
    {.}         
    {.5em}      
    {}          

\theoremstyle{boldremark}

\newtheorem{example}{Example}[section]

\newtheorem{ass}{Assumption}[section]
\newtheorem{rmk}{Remark}[section]

\numberwithin{equation}{section}
\renewenvironment{proof}{{\bfseries Proof}}{\qed}

\def\RR{\mathbb{R}}
\def\PP{\mathbb{P}}

\def\NN{\mathbb{N}}

\def\cE{{\mathcal E}}

\def\de{{\delta}}

\def\et{{\eta}}

\def\Ga{{\Gamma}}

\def\de{{\delta}}

\def\vare{{\varepsilon}}

\allowdisplaybreaks

\def\dt{\, dt}

\begin{document}

\baselineskip 15.35pt
\numberwithin{equation}{section}
\title
{{Uniform large deviation principles for SDEs under locally weak monotonicity conditions  }
}

\date{}
\author{{Jian Wang}$^1$\footnote{E-mail:20230078@hznu.edu.cn}~~~{Hao Yang}$^2$\footnote{E-mail:yanghao@hfut.edu.cn}
\\
 \small 1. School of Mathematics, Hangzhou Normal University, Hangzhou 311121, China. \\
 \small 2. School of Mathematics, Hefei University of Technology, Hefei, Anhui 230009, China.
}

\maketitle	

\newcommand\blfootnote[1]{%
\begingroup
\renewcommand\thefootnote{}\footnote{#1}%
\addtocounter{footnote}{-1}%
\endgroup
}

\begin{center}
\begin{minipage}{130mm}
{\bf Abstract:}
In this paper, we provide a criterion on uniform large deviation principles (ULDP) for stochastic differential equations under locally weak monotone conditions and Lyapunov conditions, which can be applied to stochastic systems with coefficients of polynomial growth and possible degenerate driving noises, including the stochastic Hamiltonian systems. The weak convergence method plays an important role in obtaining the ULDP. This result extends the scope of applications of the main theorem in \cite{WYZZ}.

\vspace{3mm} {\bf Keywords:}
Uniform large deviations; Lyapunov conditions; Hamiltonian
systems; weak convergence method

\vspace{3mm} {\bf AMS Subject Classification:}
60B10; 60F10; 60H10; 37A50.

\end{minipage}
\end{center}

\newpage

\renewcommand\baselinestretch{1.2}
\setlength{\baselineskip}{0.28in}

\section{Introduction}
As a motivating example, consider the family of small noise stochastic diﬀerential equations (SDEs)
\begin{equation}\label{11}
dX^\vare(t)=b(X^\vare(t))\dt+\sqrt{\vare}\sigma(X^\vare(t))dB(t),\quad t\ge 0, \quad X^\vare(0)=x \in \RR^d,
\end{equation}
where $\{B(t)\}_{t \geq 0}$ is an \emph{m}-dimensional Brownian motion on a complete probability space $(\Omega,\mathcal{F},\{{\mathcal{F}}_t\}_{t\geq 0},\mathbb{P})$   with a filtration $\{{\mathcal{F}}_t\}_{t\geq 0}$ satisfying the usual conditions, $b$ is a measurable vector field on $\RR^d$, $\sigma(\cdot): \RR^d\rightarrow \RR^d\otimes \RR^m$ is a measurable mapping, and $\vare$ is a strictly positive constant. Notice that the $X^{\vare,x}$ are indexed both by the size
of the noise $\vare$ and the initial condition $x$.

The small noise large deviation principle (LDP) for Eq. (\ref{11}) has a long history and has been studied by many authors. The general LDP was first
formulated by Varadhan \cite{vara} in 1966. Then after the pioneering work of the LDP on Markov process \cite{don&vara} and dynamical systems \cite{F-W} in the 1970s and 1980s, respectively, the LDP has attracted considerable attention as it deeply reveals the rules of the extreme events in risk management, statistical mechanics, informatics, quantum physics and so on.

For several applications, such as characterizing the exit time of $X^{\vare,x}$ from a domain, the large deviations of $X^{\vare,x}$ must be uniform with respect to the initial conditions in certain bounded subsets of the space \cite{DZ, F-W}.
The uniform large deviation principles (ULDP) are essential to the study of the dynamical behavior of the random system \cite{JWZZ}.
In this paper, our goal is to study the ULDP for Eq. \eqref{11} under locally weak monotonicity conditions and some Lyapunov conditions on the coefficients, which can be applied to stochastic systems with coefficients of polynomial growth and possible degenerate driving noises.

The classical method for LDPs relies on some discretization/approximation arguments and exponential-
type probability estimates (see, e.g. \cite{FZ}, \cite{MZ}). However, this method seems not easy to be
applied to our problem in this paper since the locally weak monotonicity conditions on the coefficients are very weak, which are even not local Lipschitz. Even if it could
be proved by the exponential equivalence method, the computation will be complicated, and
the assumptions that the coefficients satisfy some time regularity were required.

In this paper, we will adopt the weak convergence approach introduced by Budhiraja, Dupuis
and Maroulas in \cite{dupuis-1} and \cite{dupuis-1-1}. This method has proved to be
smoother and can handle stochastic partial differential equations with highly nonlinear terms . In \cite{F-W}, under either of the following circumstances:\\
(${\bf A}_1$). The coefficients  $b$ and $\sigma$ are bounded, locally Lipschitz continuous and uniformly continuous on $\mathbb{R}^d$  and there exists a positive constant $\lambda>0$ such that
$$\beta^*\sigma(x)\sigma^{*}(x)\beta\ge \lambda |\beta|^2\ {\rm for\ all}\ \beta, x\in \mathbb{R}^d;$$
(${\bf A}_2$). $d=m$, $\sigma$ is an identity matrix and the drift $b$ is globally Lipschitz continuous;
\\
Freidlin and Wentzell proved that the family of the laws of the solution processes  $\{X^{\varepsilon,x}\}$ of SDEs (\ref{11})
satisfies a ULDP with respect to the initial value $x\in \mathbb{R}^d$; see \cite [Theorem 3.1, p.135] {F-W} etc. In \cite{JWZZ}, the first
named author of this paper and his collaborators derived the ULDP for SDEs (independent of time $t$) under some locally weak monotonicity conditions and Lyapunov conditions. In this paper, we propose  more general and operational conditions on $b$ and $\sigma$, which can be applied to a wider range of examples.  As a result, Theorem \ref{thmULDP} below  covers not only the ULDP in \cite{JWZZ} but also all the examples in \cite{WYZZ}. Moreover, the stochastic Hamiltonian systems are also included, see Example 2.1 below.

Although the  weak convergence approach has  proved to be an effective way for the LDP in recent years, it is a challenging task for us to obtain the ULDP because the locally weak monotonicity conditions are very weak, which admits the cases of Hölder continuous coefﬁcients and the Gronwall type equalities are not applicable. Moreover, our setting includes not only some coefficients of polynomial growth but also degenerate driving noises. To prove the main theorem, we construct some special control functions, apply stopping time techniques and use a particularly suitable sufﬁcient condition proved in \cite{JWZZ} to verify the criteria of Budhiraja-Dupuis-Maroulas.

The layout of the present paper is as follows. In Section 2, we show assumptions and the main result.
Section 3 is devoted to establish the ULDP for Eq. \eqref{11} by using the weak convergence method.
\section{The main result}

Throughout,  we will use the following notation. Let $(\mathbb{R}^d, \langle \cdot, \cdot \rangle, |\cdot|)$ be the \emph{d}-dimensional Euclidean space with the inner product $\langle \cdot, \cdot \rangle$ which induces the norm $|\cdot|$. The norm $\| \cdot \|$ stands for the Hilbert-Schmidt norm $\| \sigma \|^2 :=  \sum_{i=1}^{d} {\sum_{j=1}^{m}{\sigma_{ij}^2}}$ for any $d \times m$-matrix $\sigma = (\sigma_{ij}) \in \mathbb{R}^d \otimes \mathbb{R}^m$. $\sigma^*$ stands for the transpose of the matrix $\sigma$.

Let $(\Omega,\mathcal{F},\{{\mathcal{F}}_t\}_{t\geq 0},\mathbb{P})$  be a complete probability space with a filtration $\{{\mathcal{F}}_t\}_{t\geq 0}$ satisfying the usual conditions and $\{B(t)\}_{t \geq 0}$  an \emph{m}-dimensional Brownian motion on this probability space. Fix $T\in(0,\infty)$. Consider the following stochastic differential equations:
\begin{equation}\label{12}
dX^\vare(t)=b(t,X^\vare(t))\dt+\sqrt{\vare}\sigma(t,X^\vare(t))dB(t),\ t\in[0,T], \quad X(0)=x \in \RR^d,
\end{equation}
where $\sigma : \mathbb{R}\otimes\mathbb{R}^d \ni (t,x) \mapsto \sigma (t,x) \in \mathbb{R}^d\otimes \mathbb{R}^m$ and $ b:  \mathbb{R}\otimes\mathbb{R}^d\ni (t,x) \mapsto b(t,x) \in \mathbb{R}^d$ are continuous. In the following, we use the notation $X^{\vare,x}$ to indicate the solution of (\ref{12}) starting from $x$.

Let us now introduce the following assumptions.
Let $f,~g,~l$ be nonnegative integrable functions on $[0,T]$.

\begin{ass}\label{ass0}
For arbitrary $R>0$,
\begin{equation*}
\int_0^T\sup_{|x|\leq R}(|b(s,x)|+\|\sigma(s,x)\|^2)ds<\infty.
\end{equation*}
\end{ass}
\begin{ass}\label{ass1}
Let $\vare_0\in(0,1)$. For arbitrary $R>0$, if $|x|\vee|y|\leq R$, there exists $L_R>0$ such that the following locally weak  monotonicity condition
\begin{equation}\label{21}
2\langle x-y,b(s,x)-b(s,y)\rangle+\|\sigma(s,x)-\sigma(s,y)\|^2\leq g(s)\eta_R(|x-y|^2),\ \forall s\in[0,T]
\end{equation}
holds for $|x-y|\leq \vare_0$, where $\eta_R:[0,1)\rightarrow \mathbb{R}_+$ is an increasing, continuous function satisfying
\begin{equation*}
\eta_R(0)=0,\quad \int_{0^+}\frac{dx}{\eta_R(x)}=+\infty.
\end{equation*}
\end{ass}

\begin{ass}\label{ass2}
There exist a Lyapunov function $V\in C^2(\mathbb R^d;\mathbb R_+)$  and $\theta>0,~\eta>0$ such that
\begin{equation}\label{As b V}
\lim_{|x|\rightarrow+\infty}V(x)=+\infty,
\end{equation}
\begin{equation}\label{22}
\langle b(s,x),\nabla V(x)\rangle+\frac{\theta}{2}\text{Trace}\big(\sigma^*(s,x)\nabla^2 V(x)\sigma(s,x)\big)+\frac{|\sigma^*(s,x)\cdot\nabla V (x)|^2}{\eta V(x)}\leq f(s)\big(1+\gamma(V(x))\big),
\end{equation}
and
\begin{equation}\label{23}
\text{Trace}\big(\sigma^{*}(s,x)\nabla^2 V(x)\sigma(s,x)\big)\geq -l(s)[M +K\gamma(V(x))].
\end{equation}
\end{ass}
\noindent Here $\gamma:[0,+\infty)\rightarrow \mathbb{R}_+$ is a continuous, increasing function satisfying
\begin{equation}\label{gamma}
\int_0^{+\infty}\frac{1}{\gamma(s)+1}ds=\infty.
\end{equation}
$\nabla V$ and $\nabla^2 V$ stand for the gradient vector and Hessian matrix of the function
 $V$, respectively, $K,\ M > 0$ are some fixed constants.
\begin{ass}\label{ass3} For any $0\leq c \leq 1$,
\begin{equation*}
\sup_{s\in [0, \varepsilon_0]}\frac{c\eta_R(s)}{\eta_R(cs)}<\infty,\quad\ \sup_{s\in [0, \infty)}\frac{c\gamma(s)}{\gamma(cs)}<\infty.
\end{equation*}
Here $\varepsilon_0$ is the constant appearing in Assumption \ref{ass1}.
\end{ass}

\begin{rmk}\label{re1}
The examples of the function {\rm  $\eta_R(s)$ in Assumption \ref{ass1} include $R s\log\frac{1}{s}$ and the examples of the function  $\gamma(s)$ in Assumption \ref{ass2}  include $s\log s+1$, etc.}
\end{rmk}

The next result gives the existence and uniqueness of the solution of SDE \eqref{12}. Its proof is is similar to that of Theorem 1.1 in \cite{wu}. The existence and uniqueness of a local solution can be established using Euler approximation.  Furthermore, one can show that the solution is global using the Lyapunov function $V$, \eqref{As b V} and \eqref{22} (cf. \cite [Theorem 3.5, p.75]{KhasminskiiB}).
\begin{prop}\label{prop1}
For any $0<\vare<1$, under Assumptions    {\rm \ref{ass0}}-{\rm \ref{ass3}}, there exists a unique  solution to Eq. \eqref{12} defined on $[0, +\infty)$.
\end{prop}
For each $h\in L^2([0,T],\mathbb R^m)$, consider the so called skeleton equation:
\begin{equation}\label{30}
dX_x^h(t)=b(t,X_x^h(t))dt+\sigma(t,X_x^h(t))h(t)dt,
\end{equation}
with the initial value $X_x^h(0)=x$.
We have the following result:
\begin{prop}\label{prop2}
Under Assumptions  {\rm\ref{ass0}}-{\rm\ref{ass3}}, there exists a unique solution to Eq. \eqref{30}.
\end{prop}
The proof of this proposition is similar to that of Proposition \ref{prop1}, so we omit it here.

\vskip 0.3cm
For any $x \in \RR^d$ and $f \in C([0,T],\mathbb{R}^d)$, we define
\begin{equation}\label{rate-0}
I_x(f)=\inf_{\left\{h \in L^2([0,T]; \mathbb{R}^m):\  f=X_x^h \right\}}\left\{\frac{1}{2}
\int_0^T|h(s)|^2ds\right\},
\end{equation}with the convention $\inf\{\emptyset\}=\infty$,  here $X_x^h\in C([0,T],\mathbb{R}^d)$ solves Eq. (\ref{30}).

We now formulate the main result on ULDP.

\begin{thm}\label{thmULDP}
For $\vare>0$, let $X^{\vare,x}$ be the solution to Eq. \eqref{12}. Suppose Assumptions {\rm \ref{ass0}}- {\rm\ref{ass3}}  are satisfied, then $I_x$ defined by \eqref{rate-0} is a rate function  on $C([0,T],\mathbb{R}^d)$ and the  family $\{I_x, x \in \RR^d\}$ of rate functions has compact
level sets on compacts. Furthermore, $\{X^{\vare,x}\}_{\vare>0}$ satisfies the ULDP on the space $C([0,T],\mathbb{R}^d)$ with the rate function $I_x$, uniformly over the initial value $x$ in bounded subsets of $\RR^d$.
\end{thm}
\begin{rmk}
If $\eta_R(z)=z$ and $\gamma(s)=s$, then Assumptions 2.1 and 2.2 in \cite{JWZZ} holds. Our hypothesis in this paper is more general than that in \cite{JWZZ}. Moreover, all the examples in \cite{WYZZ} satisfy not only the LDP, but also the ULDP.
\end{rmk}

The proof of Theorem \ref{thmULDP} will be given in Section 3 below.
Let us first look at an example of a Hamiltonian system, which is widely studied in dynamical systems, see \cite{Talay, wuliming, yong, zxc}. Obviously, (\ref{Trace}) below does not satisfy condition (2.6) in \cite{WYZZ}.

\begin{example}[\bf Hamiltonian system] Let $H : \RR^2 \rightarrow \RR$ be a
$C^2$ Hamiltonian function, which determines a Hamiltonian system by the Hamiltonian vector field $\mathcal{H}(H):= (\frac{\partial H}{\partial x_2},-\frac{\partial H}{\partial x_1})^*$.
Now we consider the following stochastic system
\begin{equation}\label{Ham}
dX(t) = \big[\mathcal{H}(H)(X(t)) -(F(H)\nabla(H))(X(t))\big]dt +\sqrt{\vare}\sigma(X(t))dB(t),\,\, \varepsilon>0,
\end{equation}
where $F: {\rm Range}(H) \rightarrow \RR$ is a $C^1$ function.

Let $H(x_1,x_2) = \frac{x_2^2}{2}+\frac{x_1^4}{4}-\frac{x_1^2}{2} \in [-\frac{1}{4}, \infty)$. Suppose that $\sigma(x)$ is a locally Lipschitz continuous function satisfying $\|\sigma(x)\|^2 \leq c_1(x_1^4+x_2^{\frac{4}{3}}+1)$, where $c_1>0$ is a fixed constant.

We consider the Lyapunov function $V = H+\frac{1}{4}$, then
$$
\nabla V(x_1,x_2)=(x_1^3-x_1,x_2)^*,\ \nabla^2V(x_1,x_2)={\rm Diag}(3x_1^2-1,1);
$$
\begin{equation}\label{Epr2}
\text{Trace}\left(\sigma^*(x)\nabla^2V(x)\sigma(x)\right)=(3x_1^2-1)(\sigma_{1,1}^2(x)+\sigma_{1,2}^2(x))+\sigma_{2,1}^2(x)+\sigma_{2,2}^2(x).
\end{equation}
Firstly, by (\ref{Epr2}), we have
\begin{eqnarray}\label{Trace}
\text{Trace}\left(\sigma^*(x)\nabla^2V(x)\sigma(x)\right)
&\ge & -\|\sigma^*(x)\|^2\notag\\
&\ge&  -c_1(x_1^4+x_2^{\frac{4}{3}}+1)\notag\\
&\ge&  -c_1(x_1^4+x_2^2+2)\notag\\
&\ge&  -4c_1-8c_1V(x),
\end{eqnarray}
that is, \eqref{23} holds.  Furthermore,
\DEQS
&&\limsup_{|x|\rightarrow +\infty}\frac{\text{Trace}\left(\sigma^*(x)\nabla^2V(x)\sigma(x)\right)}{|\nabla H(x)|^2} \\
&=& \limsup_{|x|\rightarrow +\infty}\frac{\text{Trace}\left(\sigma^*(x)\nabla^2V(x)\sigma(x)\right)}{x_1^6+x_2^2}\\
&\leq&  \limsup_{|x|\rightarrow +\infty}\frac{3x_1^2+1}{(x_1^6+x_2^2)^{\frac{1}{3}}}\times \frac{\sigma_{1,1}^2(x)+\sigma_{1,2}^2(x)}{(x_1^6+x_2^2)^{\frac{2}{3}}}+\limsup_{|x|\rightarrow +\infty} \frac{\sigma_{2,1}^2(x)+\sigma_{2,2}^2(x)}{x_1^6+x_2^2}\\
&\leq&  3c_1\limsup_{|x|\rightarrow +\infty} \frac{x_1^4+x_2^{\frac{4}{3}}+1}{(x_1^6+x_2^2)^{\frac{2}{3}}} \leq  6c_1.
\EEQS
This implies that there exists a positive constant $d_1$ such that
\begin{equation}\label{e2}
\text{Trace}\left(\sigma^*(x)\nabla^2V(x)\sigma(x)\right)\le d_1(|\nabla H(x)|^2+1)\ {\rm for\ all}\ x\in \mathbb{R}^2.
\end{equation}
Besides, the inequality
$$\limsup_{|x|\rightarrow +\infty}\frac{|\sigma^*(x)\nabla V(x)|^2}{|\nabla H(x)|^2V(x)}\le \limsup_{|x|\rightarrow +\infty}\frac{\|\sigma(x)\|^2}{V(x)}\le c_1\limsup_{|x|\rightarrow +\infty}\frac{x_1^4+x_2^{\frac{4}{3}}+1}{V(x)}\le 4c_1$$
 implies that there is a positive constant $d_2$ such that
\begin{equation}\label{e3}
\frac{|\sigma^*(x)\nabla V(x)|^2}{V(x)}\le d_2(|\nabla H(x)|^2+1)\ {\rm for\ all}\ x\in \mathbb{R}^2.
\end{equation}
Therefore, the left-hand of (\ref{22}) with $\theta =\frac{1}{d_1}, \et=2d_2$ is dominated by
\[
-F(H)(x)|\nabla H(x)|^2 + (|\nabla H(x)|^2+1).
\]
Then, we only need to choose $F$ such that $\liminf_{x \rightarrow +\infty} F(x) >1$, which means that the left-hand side of \eqref{22} is negative if $x$ is large enough. Thus, the system (\ref{Ham}) admits ULDP according to Theorem \ref{thmULDP}.
\end{example}

In the sequel, the symbol $C$ will denote a positive generic constant whose value may change from place to place.
\section{The proof}
In this section we give the proof of Theorem \ref{thmULDP}. We first recall  a sufficient condition for ULDP.
\subsection{A Sufficient Condition for ULDP}

Let $\cE = C([0,T],\mathbb{R}^d)$. $\rho(\cdot,\cdot)$ stands for the uniform metric in the space $\cE$ and $\cE_0 =\RR^d$. 
Recall that $\mathscr{K}$ is a collection of all compact subsets of $\RR^d$.

Let
\begin{equation*}
S^N :=\{h\in L^2([0,T],\mathbb R^m):|h|_{L^2([0,T],\mathbb R^m)}^2\leq N\},
\end{equation*}
and
\begin{equation*}
\tilde{S}^N :=\{\phi : \phi \ \text{is} \ \mathbb{R}^m\text{-valued} \  {\mathcal{F}}_t\text{-predictable\  process\  such\ that}\  \phi(\omega) \in S^N,\  \mathbb{P}\text{-}a.s.\}.
\end{equation*}
$S^N$ will be  endowed with the weak topology on $L^2([0,T],\mathbb R^m)$, under which $S^N$ is a compact Polish space.

For any $\vare >0$, let $\Ga^{\vare}: \cE_0 \times C([0,T],\mathbb{R}^m) \rightarrow \cE$ be a measurable mapping. Set $X^{\vare, x} := \Ga^{\vare}(x, B(\cdot))$.

The following lemma provides a sufficient condition for verifying the uniform Laplace
principle. For more detail, see Theorem 1.3 in \cite{JWZZ}.
\begin{lem}\label{suffthm}
Suppose that there exists a measurable map $\Ga^{0} : \cE_0 \times C([0,T],\mathbb{R}^m) \rightarrow \cE$ such that

{\rm (i)} for every $N < +\infty$, $x_n \rightarrow x$ and any family $\{h_n, n\in\mathbb{N}\} \subset S^N$ converging  weakly to some element $h$ as $n \rightarrow \infty$, $\Ga^0 \big( x_n, \int_0^{\cdot}{h_n}(s)ds\big)$ converges to $\Ga^0 \big(x, \int_0^{\cdot}{h(s)ds}\big)$ in the space $C([0,T],\mathbb{R}^d)$;

{\rm (ii)} for every $N < +\infty$, $\{x^{\vare}, \vare>0\} \subset \{x:  |x| \leq N \}$ and any family $\{h^{\vare}, \vare > 0\} \subset \tilde{S}^N$ and any $\delta > 0$,
\[
\lim_{\vare \rightarrow 0} \mathbb{P}\big(\rho(Y^{\vare,x^\vare},Z^{\vare,x^\vare})> \delta \big) = 0,
\]
where $Y^{\vare,x^\vare}= \Ga^{\vare}\big( x^\vare, B(\cdot) + \frac{1}{\sqrt{\vare}}\int_0^{\cdot}{h^{\vare}(s)ds}\big)$
and $Z^{\vare, x^\vare}=\Ga^0 \big(x^{\vare}, \int_0^{\cdot}{h^{\vare}(s)ds}\big)$.

Set
\begin{equation}\label{rate-0}
I_x(f)=\inf_{\left\{h \in L^2([0,T]; \mathbb{R}^m):\  f=\Ga^{0}(x,\int_0^{\cdot}h(s)ds)  \right\}}\left\{\frac{1}{2}
\int_0^T|h(s)|^2ds\right\},
\end{equation}with the convention $\inf\{\emptyset\}=\infty$. Then for all $x \in \cE_0$, $I_x$ defined by \eqref{rate-0} is a rate function on $\cE$, the family $\{I_x, x \in \cE_0\}$ of rate functions has compact level sets on compacts and $\{X^{\vare, x}\}_{\vare>0}$ satisfies a uniform Laplace principle  with the rate function $I_x$ uniformly over $\mathscr{K}$.
\end{lem}
\begin{rmk}
It is well know that a uniform Laplace principle implies the ULDP, see \cite[Proposition 14]{dupuis-2} and \cite{SA, SBD}.
\end{rmk}

\subsection{Proof of Theorem \ref{thmULDP}}
Let $S^N$ and $\tilde{S}^N$ be defined as in Subsection 3.1. According to  Proposition \ref{prop2}, there exists  a measurable mapping $\Ga^0$ from $\RR^d \times C([0,T], \mathbb{R}^m)$ to $C([0,T], \mathbb{R}^d)$ such that $X_x^{h} = \Ga^0 \big(x, \int_0^{\cdot}{h(s)ds}\big)$ for  $x \in \RR^d$ and $h\in L^2([0,T],\mathbb R^m)$.

By the Yamada-Watanabe theorem, the existence of a unique strong solution of Eq. \eqref{12} and Assumption \ref{ass1} implies that for
every $\vare>0$, there exists a measurable mapping $\Ga^{\vare}:\RR^d \times C([0,T], \mathbb{R}^m)\rightarrow C([0,T], \mathbb{R}^d)$ such that
\begin{equation*}
X^{\vare,x}=\Ga^{\vare}( x, B(\cdot) ),
\end{equation*}
and applying the Girsanov theorem, for any $N>0$ and $h^{\vare}\in \tilde{S}^N$,
\begin{equation}\label{52}
Y^{\vare,x} := \Ga^{\vare}\big(x,  B(\cdot) + \frac{1}{\sqrt{\vare}}\int_0^{\cdot}{h^{\vare}(s)ds}\big)
\end{equation}
is the solution of the following SDE
\begin{equation}\label{32}
Y^{\vare,x}(t)=x+\int_0^tb(s,Y^{\vare}(s))ds+\int_0^t\sigma(s,Y^{\vare}(s)) h^\vare(s)ds+\sqrt{\vare}\int_0^t\sigma(s,Y^{\vare}(s))dB(s).
\end{equation}

By virtue of Lemma \ref{suffthm}, to prove Theorem \ref{thmULDP}, we need to verify the conditions (i) and (ii) in Lemma \ref{suffthm} for the measurable maps $\Ga^\vare$ and $\Ga^0$.

The verification of Conditions (i) and (ii) is similar to the proof of Proposition 3.1 in \cite{WYZZ}, we here only give a sketch.

{\bf Proof of condition (i):}   Let $x_n \rightarrow x$ and $\{h_n\}_{n \in \NN} \subset S^N$ converges in the weak topology to $h$ as $n \rightarrow \infty$.

Define $\varphi(y)=\int_0^y\frac{1}{\gamma(s)+1}ds$ and $W^{h_n}_{x_n}(t)=e^{-\eta\int_0^{t}|h_n(s)|^2ds}V(X_{x_n}^{h_n}(t))$, where $\eta > 0,~\gamma,$ and $V$ are in Assumption \ref{ass2}. Apply the chain rule to get
\begin{eqnarray}
 &&\varphi(W^{h_n}_{x_n}(t))\nonumber\\
 &=& \varphi(V(x))+\int_0^{t} \varphi^\prime(W^{h_n}_{x_n}(s))\cdot e^{-\eta\int_0^s|h_n(r)|^2dr} \cdot\big[-\eta|h_n(s)|^2V(X^{h_n}_{x_n}(s)) \nonumber\\
  &&+\langle b(s,X_{x_n}^{h_n}(s)),\nabla V(X_{x_n}^{h_n}(s))\rangle
+\langle \nabla V(X_{x_n}^{h_n}(s)),\sigma(s,X^{h_n}_{x_n}(s))\cdot h_n(s)\rangle\big]ds\nonumber\\
  &\leq& \varphi(V(x))+\int_0^{t} \varphi^\prime(W^{h_n}(s))\cdot e^{-\eta\int_0^s|h_n(r)|^2dr} \cdot\big[\langle b(s,X_{x_n}^{h_n}(s)),\nabla V(X_{x_n}^{h_n}(s))\rangle  \nonumber\\
  &&+\frac{|\sigma^{*}(s,X_{x_n}^{h_n}(s)) \cdot \nabla V (X_{x_n}^{h_n}(s))|^2}{\eta V(X_{x_n}^{h_n}(s))}\big]ds\nonumber\\
  &\leq& \varphi(V(x))+\int_0^{t} \varphi^\prime(W^{h_n}_{x_n}(s))\cdot e^{-\eta\int_0^s|h_n(r)|^2dr} \cdot
 \big[ f(s)\big(1+\gamma(V(X^{h_n}_{x_n}(s)))\big)\nonumber\\
  &&+\frac{\theta}{2}l(s)\big(M+K\gamma(V(X^{h_n}_{x_n}(s))\big)\big]ds\nonumber\\
  &\leq& \varphi(V(x))+C\int_0^T(f(s)+l(s)) \cdot \frac{e^{-\eta\int_0^s|h_n(r)|^2dr}\big(1+\gamma(V(X^{h_n}_{x_n}(s)))\big)}{1+ \gamma(e^{-\eta\int_0^s|h_n(r)|^2dr}V(X^{h_n}_{x_n}(s)))}ds \nonumber\\
   &\leq& \varphi(V(x))+C\int_0^T(f(s)+l(s))ds.\label{34}
\end{eqnarray}
Assumption \ref{ass3} has been used in getting the last inequality.

 The above inequality (\ref{34}) yields
\begin{equation*}
  \sup_{n\in\mathbb{N}}\sup_{t\in[0,T]}V(X^{h_{n}}_{x_n}(t))<\infty.
\end{equation*}
By the condition (\ref{As b V}) on the function $V$, we deduce that
\begin{equation}\label{4-1}
  \sup_{n\in\mathbb{N}}\sup_{t\in[0,T]}|X^{h_n}_{x_n}(t)|\leq L
\end{equation}
for some constant $L>0$.

Next, using the Arzela-Ascoli theorem, we can show that $\{X_{x_n}^{h_n},n\in\mathbb{N}\}$ is pre-compact in the space $C([0,T],\mathbb R^d)$ and
\[
 X_{x_n}^{h_n} \rightarrow \tilde{X} \quad \text{in}~C([0,T],\mathbb{R}^d),
\]
for some $\tilde{X}$.
Then, the uniqueness of skeleton equation \eqref{30} implies $ \tilde{X} = X_x^h$. Therefore, condition (i) holds.

{\bf Proof of condition (ii):} Let $\{x^{\vare}, \vare>0\} \subset \{x:  |x| \leq N \}$ and a family $\{h^{\vare}, \vare > 0\} \subset \tilde{S}^N$, we need to prove that
$\rho(Y^{\vare,x^\vare},Z^{\vare,x^\vare}) \rightarrow 0$ in probability as $\vare \rightarrow 0$, where $Y^{\vare,x^\vare}= \Ga^{\vare}\big( x^\vare, B(\cdot) + \frac{1}{\sqrt{\vare}}\int_0^{\cdot}{h^{\vare}(s)ds}\big)$
and $Z^{\vare, x^\vare}=\Ga^0 \big(x^{\vare}, \int_0^{\cdot}{h^{\vare}(s)ds}\big)$.

Recall that $\vare_0$ is the constant appeared in Assumption \ref{ass1}.
For $R>0$, $0<p\leq \vare_0$, define stopping time
\begin{equation*}
\tau^{\vare}_R=\inf\{t\geq0:|Y^{\vare,x^\vare}(t)|\geq R\},\quad \tau^\vare_{p}=\inf\{t\geq0:|Y^{\vare,x^\vare}(t)-Z^{\vare,x^\vare}(t)|^2\geq p\}.
\end{equation*}
From the proof of \eqref{4-1} we also see that  there exists a constant $L>0$ such  that
\[
\sup_{\vare>0}\sup_{t\in[0,T]}|Z^{\vare,x^\vare}(t)|\leq L.
\]

Under Assumptions \ref{ass0}-\ref{ass3}, following a similar proof of (4.1) in Proposition 3.1 in \cite{WYZZ}, we conclude that

\begin{equation}\label{4-3}
  \lim_{\vare\rightarrow 0}\PP\big(\tau^\vare_{p}\leq \tau^{\vare}_R\wedge T\big) \leq \lim_{\vare\rightarrow 0}\PP\big(\sup_{s\leq T\wedge\tau^\vare_R\wedge\tau^\vare_{p}}|Y^{\vare, x^\vare}(s)-Z^{\vare, x^\vare}(s)|^2\geq p\big)=0.
\end{equation}

Let $\varphi(x)$ be defined as that in the proof of condition (i)  and denote
\begin{equation}
W^{h^\varepsilon}(t)=e^{-\eta\int_0^{t}|h^\varepsilon(s)|^2ds}V(Y^{\varepsilon, x^\varepsilon}(t)).
\end{equation}
 Note that $\varphi(x)$ is a concave function on the interval $[0, \infty)$ and  the second
order derivative $\varphi''(x)$ of $\varphi(x)$ in the sense of distributions is a non-positive Radon measure.

 Applying the It\^{o}--Tanaka formula to $\varphi(W^{h^\varepsilon}(t))$ gives
 \begin{eqnarray}
 &&\mathbb{E}\varphi(W^{h^\varepsilon}(T\wedge\tau^\varepsilon_R\wedge\tau^\epsilon_{p}))\nonumber\\
 &\leq& \varphi(V(x^\varepsilon))+\mathbb{E}\int_0^{T\wedge\tau^\varepsilon_R\wedge\tau^\epsilon_{p}} \varphi^\prime(W{^{h^\varepsilon}}(s))\cdot e^{-\eta\int_0^s|h^\epsilon(r)|^2dr} \cdot\big[-\eta|h^\varepsilon(s)|^2V(Y^{\varepsilon, x^\varepsilon}(s))\nonumber\\
 &&+\langle b(s,Y^{\epsilon,x^\varepsilon}(s)),\nabla V(Y^{\varepsilon, x^\varepsilon}(s))\rangle+\langle \nabla V(Y^{\varepsilon, x^\varepsilon}(s)),\sigma(s,Y^{\varepsilon, x^\varepsilon}(s))\cdot h^\epsilon(s)\rangle\nonumber\\
 &&+\epsilon\cdot\text{Trace}\big(\sigma^*(s,Y^{\varepsilon, x^\varepsilon}(s))\nabla^2 V(Y^{\varepsilon, x^\varepsilon}(s))\sigma(s,Y^\epsilon(s))\big)\big]ds\nonumber\\
  &\leq& \varphi(V(x^\varepsilon))+\mathbb{E}\int_0^{T\wedge\tau^\varepsilon_R\wedge\tau^\varepsilon_{p}} \varphi^\prime(W{^{h^\varepsilon}}(s))\cdot e^{-\eta\int_0^s|h^\varepsilon(r)|^2dr} \cdot\big[-\eta|h^\varepsilon(s)|^2V(Y^{\varepsilon, x^\varepsilon}(s))\nonumber\\
 &&+\langle b(s,Y^{\varepsilon, x^\varepsilon}(s)),\nabla V(Y^{\varepsilon, x^\varepsilon}(s))\rangle\nonumber+\frac{|\sigma^*(s,Y^{\varepsilon, x^\varepsilon}(s))\cdot \nabla V(Y^{\varepsilon, x^\varepsilon}(s))|}{\sqrt{\eta V(Y^{\varepsilon, x^\varepsilon}(s))}}\cdot \sqrt{\eta V(Y^{\varepsilon, x^\varepsilon}(s))}|h^\varepsilon(s)|\nonumber\\
 &&+\varepsilon\cdot\text{Trace}\big(\sigma^*(s,Y^{\varepsilon, x^\varepsilon}(s))\nabla^2 V(Y^{\varepsilon, x^\varepsilon}(s))\sigma(s,Y^{\varepsilon, x^\varepsilon}(s))\big)\big]ds\nonumber\\
 &\leq&\varphi(V(x^\varepsilon))+\mathbb{E}\int_0^{T\wedge\tau^\varepsilon_R\wedge\tau^\varepsilon_{p}} \varphi^\prime(W{^{h^\varepsilon}}(s))\cdot e^{-\eta\int_0^s|h^\epsilon(r)|^2dr} \cdot\big[\langle b(s,Y^{\varepsilon, x^\varepsilon}(s)),\nabla V(Y^{\varepsilon, x^\varepsilon}(s))\rangle\nonumber\\
&&+\frac{|\sigma^*(s,Y^{\varepsilon, x^\varepsilon}(s)) \cdot \nabla V(Y^{\varepsilon, x^\varepsilon}(s))|^2}{\eta V(Y^{\varepsilon, x^\varepsilon}(s))}+\varepsilon\cdot\text{Trace}\big(\sigma^*(s,Y^{\varepsilon, x^\varepsilon}(s))\nabla^2 V(Y^{\varepsilon, x^\varepsilon}(s))\sigma(s,Y^{\varepsilon, x^\varepsilon}(s))\big)\big]ds\nonumber\\
 &\leq& \varphi(V(x^\varepsilon))+\mathbb{E}\int_0^{T\wedge\tau^\varepsilon_R\wedge\tau^\varepsilon_{p}} \varphi^\prime(W^{h^\varepsilon}(s))\cdot e^{-\eta\int_0^s|h^\varepsilon(r)|^2dr} \cdot \big[ f(s)\big(1+\gamma(V(Y^{\varepsilon, x^\varepsilon}(s)))\big) \nonumber\\
 &&+ (\frac{\theta}{2} -\varepsilon)l(s)\big(M+K\gamma(V(Y^{\varepsilon, x^\varepsilon}(s)))\big)\big]ds\nonumber\\
 &\leq& \varphi(V(x^\varepsilon))+C\int_0^T (f(s)+l(s))\cdot \frac{e^{-\eta\int_0^s|h^\varepsilon(r)|^2dr}\big(1+\gamma(V(Y^{\varepsilon, x^\varepsilon}(s)))\big)}{\gamma(W^{h^\varepsilon}(s))+1} ds \nonumber\\
 &\leq& \varphi(V(x^\varepsilon))+C\int_0^T(f(s)+l(s))ds.\label{3-4}
\end{eqnarray}
The last inequality follows from the property of $\gamma$ in Assumption \ref{ass3}. By \eqref{3-4} and the definition of $\varphi$, we deduce that
\begin{equation*}
\PP\big(\tau^\varepsilon_R\leq T\wedge\tau^\varepsilon_p\big)\leq\frac{\varphi(V(x^\varepsilon))+C\int_0^T(f(s)+l(s))ds}{\int_0^{e^{-\eta N}\cdot V(R)}\frac{1}{\gamma(s)+1}ds}.
\end{equation*}
Then, by (\ref{As b V}) and (\ref{gamma}), letting $R\rightarrow\infty$, we obtain
\begin{equation}\label{4-4}
\lim_{R\rightarrow\infty}\sup_{\vare\in(0,1)}\PP\big(\tau^\vare_R \leq T\wedge\tau^\vare_{p}\big)=0.
\end{equation}
Finally, for arbitrary $\delta>0$, we have
\begin{eqnarray*}
&&\PP\big(\sup_{0\leq s\leq T}|Y^{\vare, x^\vare}(s)-Z^{\vare, x^\vare}(s)|\geq\de\big)\\
&=&\PP\big(\sup_{0\leq s\leq T}|Y^{\vare, x^\vare}(s)-Z^{\vare, x^\vare}(s)|\geq\de,\tau^\vare_R\wedge\tau^\vare_{\de^2}>T\big) \\
&&+\PP\big(\sup_{0\leq s\leq T}|Y^{\vare, x^\vare}(s)-Z^{\vare, x^\vare}(s)|\geq\de,\tau^\vare_R\leq T\wedge\tau^\vare_{\de^2}\big)\\
&&+\PP\big(\sup_{0\leq s\leq T}|Y^{\vare, x^\vare}(s)-Z^{\vare, x^\vare}(s)|\geq\de,\tau^\vare_{\de^2}\leq \tau^\vare_R\wedge T\big)\\
&\leq&\PP\big(\sup_{0\leq s\leq T\wedge\tau^\vare_R\wedge\tau^\vare_{\de^2}}|Y^{\vare, x^\vare}(s)-Z^{\vare, x^\vare}(s)|^2\geq \de^2\big)+\PP\big(\tau^\vare_R\leq T\wedge\tau^\vare_{\de^2}\big)+\PP\big(\tau^\vare_{\de^2}\leq \tau^\vare_R\wedge T\big).
\end{eqnarray*}
\eqref{4-3} (with $p=\de^2$) implies that
\begin{eqnarray*}
\lim_{\vare\rightarrow 0}\PP\big(\sup_{0\leq s\leq T}|Y^{\vare, x^\vare}(s)-Z^{\vare, x^\vare}(s)|\geq\de\big)
\leq
\sup_{\vare\in(0,1)}\PP\big(\tau^\vare_R\leq T\wedge\tau^\vare_{\de^2}\big).
\end{eqnarray*}
Let $R\rightarrow \infty$ and \eqref{4-4} to get
\begin{equation*}
\lim_{\vare\rightarrow0}\PP\big(\sup_{0\leq s\leq T}|Y^{\vare, x^\vare}(s)-Z^{\vare, x^\vare}(s)|\geq\de\big)=0.
\end{equation*}
Therefore, condition (ii) holds.

\vskip 0.4cm
\noindent{\bf Acknowledgements.}
This work is partially supported by  the
National Natural Science Foundation of China (No. 12401175), the Fundamental Research Funds for the Central Universities (Nos. JZ2023HGQA0118, JZ2023HGTA0171).

\vskip 0.4cm
\noindent{\bf Declarations} The authors declare that they have no conflict of interest.

\end{document}